\documentclass[a4paper,titlepage,10pt]{article}
\usepackage[english]{babel}
\usepackage[latin1]{inputenc}
\usepackage[T1]{fontenc}
\usepackage{amsmath}
\usepackage{amscd}
\usepackage{amsfonts} 
\usepackage{amssymb}
\usepackage{bbm}
\usepackage{yfonts}
\usepackage{stmaryrd}
\usepackage{amsthm}
\usepackage{cases}
\usepackage{wasysym}
\usepackage{mathrsfs}
\usepackage[all]{xy}
\usepackage{xcolor}
\usepackage{titlesec}
\usepackage{lipsum}
\titleformat{\section}{\normalfont\bf\center}{\thesection}{1em}{}
\title{}
\author{}
\date{}
\theoremstyle{plain}
\newtheorem{teorema}{Theorem}[section]
\newtheorem{corollario}[teorema]{Corollary}
\newtheorem{proposizione}[teorema]{Proposition}
\newtheorem{lemma}[teorema]{Lemma}
\theoremstyle{definition}
\newtheorem{definizione}[teorema]{Definition}
\theoremstyle{remark}
\newtheorem{osservazione}[teorema]{Remark}
\newtheorem{esempio}{Example}[section]

\newcommand{\Ad}{\text{Ad}}

\renewcommand{\d}{\text{d}}

\renewcommand{\O}{\text{\bf O}}
\newcommand{\SO}{\text{\bf SO}}

\newcommand{\U}{\text{\bf U}}
\newcommand{\SU}{\text{\bf SU}}

\newcommand{\codim}{\text{codim}}

\newcommand{\chm}{\text{chm}}

\newcommand{\ac}{\text{ac}}

\newcommand\blfootnote[1]{%
  \begingroup
  \renewcommand\thefootnote{}\footnote{#1}%
  \addtocounter{footnote}{-1}%
  \endgroup
}

\newcommand{\beq}{\begin{equation}}
\newcommand{\eeq}{\end{equation}}
\newcommand{\beqn}{\begin{equation*}}
\newcommand{\eeqn}{\end{equation*}}
\newcommand{\bpr}{\begin{proof}}
\newcommand{\epr}{\end{proof}}
\newcommand{\ben}{\begin{enumerate}}
\newcommand{\een}{\end{enumerate}}
\newcommand{\bit}{\begin{itemize}}
\newcommand{\eit}{\end{itemize}}

\newcommand{\ga}{\mathfrak{a}}

\newcommand{\gs}{\mathfrak{s}}
\newcommand{\gt}{\mathfrak{t}}

\newcommand{\bC}{\mathbb{C}}

\newcommand{\bR}{\mathbb{R}}

\newcommand{\bT}{\mathbb{T}}

\newcommand{\bZ}{\mathbb{Z}}

\newcommand{\bdefin}{\begin{definizione}}
\newcommand{\edefin}{\end{definizione}}
\newcommand{\bteor}{\begin{teorema}}
\newcommand{\eteor}{\end{teorema}}
\newcommand{\blem}{\begin{lemma}}
\newcommand{\elem}{\end{lemma}}
\newcommand{\bcor}{\begin{corollario}}
\newcommand{\ecor}{\end{corollario}}
\newcommand{\boss}{\begin{osservazione}}
\newcommand{\eoss}{\end{osservazione}}
\newcommand{\bprop}{\begin{proposizione}}
\newcommand{\eprop}{\end{proposizione}}
\newcommand{\bes}{\begin{esempio}}
\newcommand{\ees}{\end{esempio}}
\newcommand{\bdim}{\begin{proof}}
\newcommand{\edim}{\end{proof}}

\begin{document}
\renewcommand{\theenumi}{\em\roman{enumi}}
\begin{center}
{\bf \Large \scshape Representations of abstract copolarity one and two}

$\,$

{\scshape Francesco Panelli and Carolin Pomrehn}
\end{center}

$\,$

$\,$

\begin{changemargin}{0.5cm}{0.5cm} \begin{small}
\noindent {\scshape Abstract}. We study the cohomogeneity of possibly reducible representations of copolarity $1$ and $2$, generalizing results obtained by Claudio Gorodski, Carlos Olmos, Ruy Tojeiro and by Claudio Gorodski, Alexander Lytchak.
\end{small}\end{changemargin} 

$\,$

\section{Introduction}\blfootnote{2010 {\em Mathematics Subject Classification}. Primary 57S15, 20G05.}\blfootnote{The first author was supported by GNSAGA of INdAM -- Roma (Italy).}
Consider an orthogonal representation $\rho:H\rightarrow\O(W)$, written $\rho=(H,W)$ for short, of a compact Lie group $H$ on a real vector space $W$. A subspace $\Sigma$ of $W$ is called a {\em generalized section}, or, more precisely, a $k$-{\em section}, if it intersects all orbits, and contains the normal space to principal orbits it meets at each intersection point with codimension $k$ (cf. \cite{GOT}). The minimum of all such $k$, c$(\rho)$, is called the {\em copolarity} of $\rho$. The copolarity of a representation has been related to its cohomogeneity in some particular cases: for instance, it is proved that an irreducible representation of copolarity $1$ has cohomogeneity $3$ (cf. \cite[Theor. 1.1]{GOT}). In \cite[Cor. 1.6]{GL} this result was generalized showing that an irreducible representation of copolarity $k$, $1\leq k\leq 6$, has cohomogeneity $k+2$. 

In this paper we analyze the question whether these results continue to hold for non-necessarily irreducible representations in the special case $k=1,2$, giving another partial answer to the question stated by Claudio Gorodski and Alexander Lytchak in \cite[Question 1.15]{GL}, asking for general relationships between copolarty and cohomogeneity of representations. In the reducible case, there is a standard method to construct counterexamples (cf. \cite{GOT}). Indeed consider a product representation $\rho:=(H_1\times H_2,W_1\oplus W_2)$ (here $H_1\times H_2$ acts componentwisely) where $(H_2,W_2)$ is polar. Then $\text{c}(\rho)=\text{c}(H_2,W_2)$, but the cohomogeneity of $\rho$ is grater than that of $(H_1,W_1)$. In this way we have a family of representations with constant copolarity, but arbitrarily large cohomogeneity. Therefore, the question is interesting only for a smaller class of representations, namely those which are {\em indecomposable}; roughly speaking, they are representations that cannot be written as a product of subrepresentations (cf. Definition \ref{decomposable_reps} below). In the case c$(\rho)=1$ we show that the counterexamples described above are the only ones that can occour.

\bteor\label{copolarity_1}
Let $\rho:H\rightarrow\O(W)$ be an indecomposable representation of a connected, compact Lie group $H$ which is not the $1$-dimensional torus $\bT^1$. If $\rho$ has copolarity $1$, then it has cohomogeneity $3$.
\eteor

The case c$(\rho)=2$ is slightly more complicated, and in fact there are some exceptions which will be completely described. Actually our result will be proved in a more general setting, and before stating it we need to recall some definitions from \cite[Section 1]{GL}.

\bdefin
Let $\rho:H\rightarrow\O(W)$ and $\tau:G\rightarrow\O(V)$ be representations of two compact Lie groups $G$, $H$.
\ben
\item We say that $\rho$ and $\tau$ are {\em quotient-equivalent} if their orbit spaces $W/H$ and $V/G$ are isometric. If moreover $\dim G<\dim H$, we say that $\tau$ is a {\em reduction} of $\rho$.
\item We say that $\tau$ is {\em reduced} if $\dim G$ is minimal in the quotient-equivalence class of $\tau$. In this case $\ac(\tau):=\dim G$ is called the {\em abstract copolarity} of $\tau$ (or of any representation quotient-equivalent to $\tau$).
\item A reduction of $\rho$ which is reduced is called a {\em minimal reduction} of $\rho$.
\een
\edefin

The simplest example of quotient-equivalence is {\em orbit-equivalence}: two representations $\rho_1$, $\rho_2$ are called {\em orbit-equivalent} if they have the same orbits up to an isometry between the representation spaces; in this case we write $\rho_1\simeq_{o.e.}\rho_2$. 

For later use, we shall now briefly describe two simple ways to obtain a reduction from a given representation $\rho:H\rightarrow\O(W)$.

First consider the principal isotropy group $K$ of $\rho$, and let $W^K$ denote the subspace of all fixed points of $K$. Clearly the normalizer $N_H(K)$ acts on $W^K$. The inclusion map $W^K\rightarrow W$ induces then an isometry $W^K/\bar{N}\rightarrow W/H$, where $\bar{N}:=N_H(K)/K$ (cf. \cite{GT, LR, Str}). If $K$ is not trivial, then $(\bar{N}, W^K)$ is a reduction of $(H,W)$, which is called {\em Luna-Richardson-Straume} reduction, LRS-reduction in the sequel.

Let now $\Sigma$ be a $k$-section for $\rho$, where $k$ is minimal (so that $k=\text{c}(\rho)$). The group $N_H(\Sigma)$ of all elements $h\in H$ preserving $\Sigma$ acts on $\Sigma$ with a kernel which we denote by $Z_H(\Sigma)$. Then the inclusion map $\Sigma\rightarrow W$ yields an isometry $\Sigma/\bar{N}_\Sigma\rightarrow W/H$, where $\bar{N}_\Sigma:=N_H(\Sigma)/Z_H(\Sigma)$ (cf. \cite{GOT}). If $\Sigma$ is a proper subspace of $W$ (in this case, following \cite{GOT}, we say that $\rho$ has {\em non-trivial copolarity}), then $(\bar{N},\Sigma)$ is a reduction of $(H,W)$.
Notice that abstract copolarity is bounded above by copolarity; however it is not known whether these two invariants coincide in general (cf. \cite[p. 69]{GL}).

We may now state our second main result:

\bteor\label{cohomogeneity_abstract_copolarity_2}
Let $\rho:H\rightarrow\text{\em\bf O}(W)$ be a non-reduced, indecomposable representation of a compact connected Lie group $H$ of abstract copolarity $2$. Then either 
\ben
\item $\rho$ has cohomogeneity $4$, or
\item $\rho=\rho_1\oplus\rho_2$ where 
\ben
\item $\rho_1:H\rightarrow\O(W_1)$ is orbit-equivalent to the isotropy representation of a rank $2$ real grassmannian,
\item $\rho_2:H\rightarrow\O(W_1^\perp)$ is orbit-equivalent to a non-polar $\bT^1$-representation without non-trivial fixed points.
\een\een

\noindent Conversely, let $\rho_1$ be the isotropy representation of a rank $2$-real grasmannian, $\rho_2$ be a non-polar $\bT^1$-representation without non-trivial fixed points and set $\rho:=\rho_1\oplus\rho_2$. Then $\rho$ is indecomposable, has cohomogeneity $\neq 4$, and both its copolarity and abstract copolarity are equal to $2$.
\eteor

We remark here that Theorem \ref{copolarity_1} may be expressed in terms of abstract copolarity as well, and in fact it will be proved in such a setting; we refer to Section \ref{Proof_of_the_main_results}, Theorem \ref{cohomogeneity_abstract_copolarity_1}, for its precise statement. From the latter result we shall be able to derive the following:

\bcor\label{c1_iff_ac1}
A representation $\rho:H\rightarrow\O(W)$ of a connected compact Lie group $H$ has abstract copolarity $1$ if and only if it has copolarity $1$.
\ecor

This means in particular that Theorems \ref{copolarity_1} and \ref{cohomogeneity_abstract_copolarity_1} are equivalent. We don't know whether the analogous of Corollary \ref{c1_iff_ac1} when copolarity and abstract copolarity are equal to $2$ holds; however note that a representation of a connected compact Lie group which has copolarity $2$ must have abstract copolarity $2$. 

Another motivation for Theorems \ref{copolarity_1}, \ref{cohomogeneity_abstract_copolarity_2} comes from Theorem 1.7 in \cite{GL}. Indeed, an easy consequence of such a result asserts that if $\rho=(H,W)$ is an irreducible representation of a connected compact group $H$ admitting a minimal reduction $\tau=(G,V)$ where $G$ is a finite extension of a $k$-dimensional torus, $k\geq 1$, then the cohomogeneity of $\rho$ equals $k+2$.

As a technical tool we shall need the following observation, which, we think, is interesting on its own, and might be useful in other contexts:

\bprop\label{identity_rep_indec}
Let $\rho:H\rightarrow\O(W)$ be a non-polar indecomposable representation of a connected, compact Lie group $H$. If $G\rightarrow\O(V)$ is quotient-equivalent to $\rho$, then the induced representation $G^\circ\rightarrow\O(V)$ is indecomposable.
\eprop

The paper is organized as follows. In Section \ref{Decomposability_of_representations} we formally define decomposable representations and give a useful criterion for decomposability of toric representations. In Section \ref{Nice_Involutions} we prove Proposition \ref{identity_rep_indec} and recall the notion of {\em nice involution} from \cite{GL}, adapting it to our particular setting. Finally in Section \ref{Proof_of_the_main_results} we prove our main results and Corollary \ref{c1_iff_ac1}.

$\,$

\noindent{\bf Acknowledgements.} The first author is very grateful to Prof. A. Lytchak for kindly introducing him into the topic of this paper; he wishes also to thank Prof. C. Gorodski for very useful conversations and many valuable remarks and Prof. F. Podest\`a for several suggestions and his constant support. Part of this work was completed while the first author was visiting the University of K\"oln, which he wishes to thank for hospitality.

The second author would like to thank Prof. G. Thorbergsson, Prof. A. Lytchak and Dr. S. Wiesendorf for their support and many helpful discussions.

\section{Decomposability of representations}\label{Decomposability_of_representations}
In this section we rigorously define the notion of {\em decomposability} for representations and develop some properties and consequences of the definition that will be needed later. 

\bdefin\label{decomposable_reps}
Let $\rho=(H,W)$ be a representation of a compact Lie group $H$. We say that $\rho$ is {\em decomposable} if the metric space $W/H$ splits as a product
$$
W/H=X_1\times X_2
$$
for suitable metric spaces $X_1$, $X_2$ neither of which is isometric to a single point.
\edefin

Note that the notion of decomposability has already been used in the more general setting of singular Riemannian foliations by Alexander Lytchak in \cite{Lyt14} and Marco Radeschi in \cite{Rad14}.

\blem\label{dec_1}
Let $(H,W)$ be a decomposable representation of the compact Lie group $H$, and set $W/H=X_1\times X_2$. Then there exists an $H$-invariant splitting $W=W_1\oplus W_2$ of $W$ such that $X_i$ is isometric to the quotient space of the induced representation $\rho_i:=(H,W_i)$, $i=1,2$. In particular $(H,W)$ is orbit-equivalent to the product representation $(H_1\times H_2,W_1\oplus W_2)$, where $H_i:=\rho_i(H)$.
\elem
\bdim
Let $SW$ denote the unit sphere in $W$, and set $S:=SW/H$, $X:=W/H$. Notice that $S$ is the unit sphere in $X$, which in turn coincides with the cone $CS$ over the metric space $S$. So we have a splitting of the cone $CS$ as a product $X_1\times X_2$. Let $O$ be the vertex of $CS$, and use the same letter to denote both the projections of $O$ in $X_1$ and $X_2$. Identify then $X_i$ with the $X_i$-fibre through $O$, and set $Y_i:=X_i\cap S$, $i=1,2$. It is easily seen that $X_i$ is isometric to the cone $CY_i$ over $Y_i$, so $S$ is isometric to the spherical join $Y_1\ast Y_2$ of $Y_1$ and $Y_2$ (cf. \cite[p. 63]{BH} and note that diam$(S)\leq\pi$). Since, for $i,j=1,2$, $i\neq j$, $Y_i$ is the set of all points in $S$ having dinstance $\pi/2$ from each point in $Y_j$, we deduce, using the discussion in \cite[p. 76]{GL}, that we may find two $H$-invariant complementary subspaces $W_1$, $W_2$ of $W$ such that $Y_i=SW_i/H$, i.e. $X_i=W_i/H$.

Now we show that the product representation $(H_1\times H_2,W_1\oplus W_2)$ and $(H,W)$ have the same orbits. Indeed, if $w \in W$, it is clear that $H\cdot w\subseteq (H_1\times H_2)\cdot w$. On the other hand, the representations above are quotient-equivalent, hence there exists $w'\in W$ such that $(H_1\times H_2)\cdot w'=H\cdot w$. Thus $(H_1\times H_2)\cdot w'=H\cdot w\subseteq (H_1\times H_2)\cdot w$ which yields $H\cdot w = (H_1\times H_2)\cdot w$, as claimed.
\edim

By the main theorem in \cite{FL08}, given a representation $(H,W)$ the quotient metric space $W/H$ splits as a product of metric spaces
\beq\label{de_rham}
W/H=X_0\times X_1\times \cdots\times X_n,
\eeq
where $X_0$ is euclidean (possibly trivial), and where each $X_i$, $i\geq 1$, is a non-trivial non-euclidean metric space which is irreducible, in the sense that it cannot be further decomposed as a non-trivial product of metric spaces. Moreover decomposition (\ref{de_rham}) is unique up to a permutiation of the isometric factors. Lemma \ref{dec_1} yields then:

\bprop\label{deRham_dec_quot}
Let $(H,W)$ be a representation of a compact Lie group $H$. Then there exists an $H$-invariant splitting $W=W_0\oplus W_1\oplus\cdots\oplus W_n$ of $W$ such that:
\ben
\item the induced representation $\rho_0:=(H,W_0)$ is trivial and $\dim W_0\geq 0$;
\item the induced representation $\rho_i:=(H,W_i)$ is indecomposable, non-trivial and $\dim W_i\geq 1$ for  $i=1,\dots,n$;
\item $(H,W)$ is orbit-equivalent to the product representation $$(H_0\times H_1\times\cdots\times H_n,W_0\oplus W_1\oplus\cdots\oplus W_n),$$ where $H_0:=\{1\}$, $H_i:=\rho_i(H)$, $i=1,\dots,n$.
\een
\eprop

\boss\label{ac_of_a_product}
The abstract copolarity of a product representation $\rho$ of the form $(H_1\times H_2,W_1\oplus W_2)$ equals the sum of the abstract copolarities of $\rho_i:=(H_i,W_i)$, $i=1,2$. Indeed, let $(G,V)$ be a minimal reduction of $\rho$, and let $(G_i,V_i)$ be a minimal reduction of $\rho_i$. By Lemma \ref{dec_1} $(G,V)$ is orbit-equivalent to a product representation $(G_1'\times G_2',V_1'\oplus V_2')$, where $V=V_1'\oplus V_2'$. Clearly we have
$$
\dim G\leq \dim G_1+\dim G_2\leq \dim G_1'+\dim G_2'.
$$
The first inequality means that $\ac(\rho)\leq\ac(\rho_1)+\ac(\rho_2)$. If the two sides were not equal, than $\dim G<\dim G_1'+\dim G_2'$, showing that $(G,V)$ is a minimal reduction of $(G_1'\times G_2',V_1'\oplus V_2')$; so $\dim V<\dim V_1'+\dim V_2'$, a contradiction.
\eoss

We shall use now Lemma \ref{dec_1} to characterize decomposability of representations of the torus $\bT^k$. This result will be needed several times later on.

It is well-known that (finite-dimensional) irreducible complex representations of $\bT^k$ are in 1-1 correspondence with Lie groups homomorphisms $\bT^k\rightarrow S^1$. Moreover, given any complex representation $(\bT^k,U)$ and any such homomorphism $\phi$, we denote by $U_\phi$ the {\em isotypical component} of $\phi$, i.e. the generalized eigenspace of $U$ with generalized eigenvalue $\phi$. If $U_\phi\neq \{O\}$, we say that $\phi$ is a {\em character} of the representation. Clearly $U$ can be uniquely decomposed as a direct sum of isotypical components. 

In this paper we are interested only in real representations $\rho=(\bT^k, V)$, therefore we apply the theory to the complexified representation $(\bT^k,V^\bC)$. Let $V_0$ be the set of fixed points of $\rho$; then $V$ can be uniquely decomposed as
$
V=V_0\oplus V_1\oplus\cdots\oplus V_n,
$
where, $V_i$, $i\geq 1$, are the {\em (real) isotypical components} of $V$, i.e. each $V_i$ has the form $(V^\bC_\phi\oplus V^\bC_{\bar{\phi}})\cap V$ for a suitable non-trivial character $\phi$ of $\rho$.

We identify the Lie algebra of $S^1$ with the real line, and associate to any character $\phi$ of $\rho$ (and so to any isotypical component of $\rho$) a $1$-form $\theta_\phi$ on the Lie algebra $\gt^k$ of $\bT^k$, namely $\theta_\phi:=\d\phi_e$. Such $1$-forms are called {\em weights} of $\rho$.

\bprop\label{two_subsets}
Let $\rho=(\bT^k,V)$ be a representation of the torus $\bT^k$, and let $\Theta\subseteq(\gt^k)^\ast$ denote the set of all weights of $\rho$. Assume that there exist two non-empty subsets $\Theta_1$, $\Theta_2\subseteq\Theta$ so that $\Theta_1\cup\Theta_2=\Theta$ and $\langle \Theta_1\rangle\cap\langle \Theta_2\rangle=\{O\}$. Then $\rho$ is decomposable.
\eprop
\bdim
If $\Theta_i$, $i=1,2$, are as in the statement, we denote by $V_i$ the sum of all isotypical components of $\rho$ corresponding to weights in $\Theta_i$, and by $\ga_i\subseteq\gt^k$ the annihilator of $\langle\Theta_i\rangle$. From $\langle \Theta_1\rangle\cap\langle \Theta_2\rangle=\{O\}$ we get $\ga_1+\ga_2=\gt^k$ (direct sum if and only if $\rho$ has discrete kernel), while from $\Theta_1\cup\Theta_2=\Theta$ we get $V_1\oplus V_2=V$. Let $T_i$ be the subtorus of $\bT^k$ with Lie algebra $\ga_i$. The exponential map of $\bT^k$ is surjective, hence any element $t\in \bT^k$ can be written as $t=t_1t_2$ for suitable $t_i\in T_i$. Since the induced representation $(T_i,V_i)$ is trivial, we deduce that $\rho$ and the product representation $(T_2\times T_1,V_1\oplus V_2)$ have the same orbits, so $\rho$ is decomposable.
\edim

\boss
Using notation of Proposition \ref{two_subsets}, we see that $\rho$ has discrete kernel if and only if $\langle\Theta\rangle=(\gt^k)^\ast$. Indeed, an element $x\in\gt^k$ belongs to the Lie algebra of $\ker\rho$ if and only if it is annihilated by all weights of $\rho$.
\eoss

In order to state the main Corollary of Proposition \ref{two_subsets} we need to introduce some terminology. Precisely, if $\rho=(\bT^k,V)$ is a representation of the torus $\bT^k$, we say that a line $\gs\subseteq(\gt^k)^\ast$ is {\em induced} by $\rho$ if $\rho$ has a character $\phi$ such that $\langle\theta_\phi\rangle=\gs$.

\bcor\label{k_lines}
Let $\rho=(\bT^k,V)$ be a faithful representation of the torus $\bT^k$, $k\geq 2$. If $\rho$ is indecomposable, then it induces at least $k+1$ lines in $(\gt^k)^\ast$.
\ecor
\bdim
Let $\gs_1,\dots$, $\gs_\ell$, $\ell\leq k$, be the lines induced by $\rho$ in $(\gt^k)^\ast$, and let $\Theta$ be the set of all weights of $\rho$. Clearly $\langle\gs_1,\dots,\gs_\ell\rangle=\langle\Theta\rangle=(\gt^k)^\ast$, where the latter equality holds because $\rho$ is faithful, hence $\ell=k\geq 2$, and the $\gs_i$'s are linearly independent. Now define $\Theta_1$ to be the subset of $\Theta$ consisting of all weights $\theta$ such that $\langle\theta\rangle=\gs_1$, and similarly define $\Theta_2$ to be the subset of $\Theta$ consisting of all weights $\theta$ such that $\langle\theta\rangle=\gs_j$, for some $j\geq 2$. Then $\Theta_1$ and $\Theta_2$ satisfy the hypotesis of Proposition \ref{two_subsets} therefore $\rho$ is decomposable.
\edim

Especially for toric representations, the absence of proper generalized sections can be read of the quotient. In order to state precisely this fact, recall that, given an orthogonal representation $(H,W)$ of a compact Lie group $H$, the {\em stratum} through a point $p\in W$, Str$(p)$, is, by definition, the connected component through $p$ of the set of points whose isotropy groups are conjugated to the isotropy at $p$. We have
\beq\label{dim_stratum}
\dim\text{Str}(p)=\dim W^{H_p}+\dim H-\dim N_H(H_p),
\eeq 
where $H_p$ denotes the isotropy at $p$, while $W^{H_p}$ denotes the set of fixed points of $H_p$. The {\em boundary} (in the sense of Alexandrov) of $W/H$, $\partial(W/H)$, is defined as the closure of the union of all strata in $W/H$ that have codimension $1$.

\bprop\label{boundary_and_copolarity_toric_actions}
Let $\rho=(\bT^k,V)$ be a representation of the $k$-dimensional torus $\bT^k$. Then $\partial(V/\bT^k)= \varnothing$ if and only if $\rho$ has trivial copolarity.
\eprop
\bdim
We may assume without loss of generality that $\rho$ is faithful and has no non-trivial fixed points. Since $\bT^k$ is abelian, this in particular implies that $\rho$ has trivial principal isotropy group.

If $\partial(V/\bT^k)= \varnothing$, $\rho$ is reduced by \cite[Prop. 1.1]{GL}, so it has trivial copolarity.

Conversely, assume $\partial(V/\bT^k)\neq \varnothing$, and let $p\in V$ be a point projecting onto a stratum of codimension $1$ in the quotient. The non-trivial part of the slice representation at $p$ has trivial principal isotropy group and cohomogeneity $1$, therefore $\bT^k_p$ is a sphere $S^a$. Since $\bT^k$ is abelian we get $a\in\{0,1\}$ and, by (\ref{dim_stratum}), $\dim V^{\bT^k_p}=\dim\text{Str}(p)=\dim V-a-1$, where $V^{\bT^k_p}$ denotes the fixed point set of $\bT^k_p$.

Assume first $a=0$, and let $\omega$ be the generator of $\bT^k_p\simeq S^0\simeq\bZ_2$. Then  $V^\omega$ has codimension $1$ and $\omega$ is a reflection. On the other hand $\omega\in S^1\subseteq\SO(V)$, a contradiction.

Then $a=1$, and we have a $\bT^k$-invariant decomposition
$
V=V^{\bT^k_p}\oplus\bar{V},
$
where $\dim\bar{V}=2$; we shall write accordingly $\rho=\tilde{\rho}\oplus\bar{\rho}$. Notice that $\bar{\rho}$ is irreducible since $\rho$ has no non-trivial fixed points. Let $\tilde{\Theta}$ be the set of weights corresponding to isotypical components in $V^{\bT^k_p}$, and denote by $\bar{\theta}$ the weight of $\bar{\rho}$. Clearly any weight of $\rho$ either belongs to $\tilde{\Theta}$, or it is equal to $\bar{\theta}$; moreover all weights in $\tilde{\Theta}$ vanish on the $1$-dimensional Lie algebra of $\bT^k_p$, $\gt^k_p$. Faithfulness of $\rho$ implies that $\bar{\theta}$ cannot vanish on $\gt^k_p$, so $\langle\tilde{\Theta}\rangle\cap\langle\bar{\theta}\rangle=\{O\}$. By (the proof of) Proposition \ref{two_subsets} $\rho$ is orbit-equivalent to the product representation $(\tilde{\rho}(\bT^k)\times S^1,V^{\bT^k_p}\oplus \bar{V})$; since $(S^1,\bar{V})$ is polar $\rho$ cannot have trivial copolarity. 
\edim

\section{Nice involutions}\label{Nice_Involutions}
In this paper an important role is played by representations $(G,V)$ of a compact (possibly disconnected) Lie group $G$ such that $\Gamma:=G/G^\circ$ acts on $V/G^\circ$ as a reflection group, i.e. $\Gamma$ is generated by elements which act on $V/G^\circ$ as reflections. Notice that the representation $(G,V)$ appearing in the statement of Proposition \ref{identity_rep_indec} is of this kind, since it is quotient-equivalent to a representation of a connected group (cf. \cite[Prop. 1.2]{GL}).

\bdim[Proof of Proposition \ref{identity_rep_indec}]
By hypotesis and what we observed above, $(G,V)$ is non-polar, indecomposable and the quotient group $\Gamma:=G/G^\circ$ is generated by reflections. Decompose the quotient $V/G^\circ$ of the induced representation $(G^\circ,V)$ as in (\ref{de_rham}):
\beq\label{a}
V/G^\circ= X_0\times X_1\times\cdots \times X_k.
\eeq
We may assume that $X_1,\dots$, $X_m$ are the non-flat factors appearing in (\ref{a}), and set $X_{\text{flat}}:=X_0\times X_{m+1}\times\cdots\times X_k$. Since a representation is polar if and only if the principal part of its orbit space is flat (cf. \cite{Ale, HLO}), and since non-polarity of $(G,V)$ is equivalent to non-polarity of $(G^\circ,V)$, we deduce that $m\geq 1$. Now, any reflection in $\Gamma$ either preserves all $X_i$, being a reflection on one of them and fixing pointwisely the others, or interchanges $X_j$, $X_m$, $j\neq m$, and fixes pointwisely the others. In the latter case $X_j$, $X_m$ have dimension $1$, so they are flat; therefore
$$
\frac{V}{G}=\frac{X_{\text{flat}}}{\Gamma}\times\frac{X_{1}}{\Gamma}\times\cdots\times\frac{X_m}{\Gamma}.
$$
Since $(G,V)$ is indecomposable and $m\geq 1$ we deduce that $X_{\text{flat}}$ is trivial and that $m=1$; this finishes the proof.
\edim

Among all representations $\rho=(G,V)$ where the quotient group $\Gamma:=G/G^\circ$ is generated by reflections, those with trivial principal isotropy group will be of particular importance to us; indeed, in this case, the discussion in \cite[Section 4.3]{GL} implies that there exists a set $\Omega\subseteq \Gamma$ of reflections generating $\Gamma$ all of whose elements $\omega'$ have a lift $\omega$ which is an involution and satisfies the following dimension formula:
\beq\label{nice_involutions}
\dim V^\omega+\dim G-\dim Z_G(\omega)=\dim V-1,
\eeq
where $V^\omega$ is the fixed-point set of $\omega$, and $Z_G(\omega)$ is the centralizer of $\omega$ in $G$. Following \cite{GL}, we shall call involutions of this kind {\em nice involutions}. 

From now on $(G,V)$ will be a minimal reduction of a non-reduced indecomposable representation $(H,W)$ with $H$ connected. Note that, in this case, the orbit space $V/G=W/H$ has non-empty boundary (cf. \cite[Prop. 1.1]{GL}). Of course, moreover, $(G,V)$ has trivial copolarity and trivial principal isotropy group. We shall denote by $\Gamma$ the finite group $G/G^\circ$.

In the following sections $G$ will be assumed to have dimension $1$ or $2$, so $G^\circ$ will be a $1$ or a $2$-dimensional torus; however in the remaining part of this Section we shall prove some preliminary results which don't require the above restriction on $\dim G$, so we shall only suppose that $G$ is a finite extension of a $k$-dimensional torus $\bT^k$, $k\geq 1$.

We begin collecting a few simple remarks in the following:

\blem\label{easy_remarks}
Under the above assumptions, the group $G$ is disconnected, and contains a set of nice involutions whose projection in $\Gamma$ is a set of generators. Moreover the induced representation $(\bT^k,V)$ of the identity component of $G$ is indecomposable and reduced, and its orbit space $V/\bT^k$ has empty boundary.
\elem
\bdim
First we observe that $G$ cannot be connected by Proposition \ref{boundary_and_copolarity_toric_actions}. Since $G$ has trivial principal isotropy group and $\Gamma$ acts on $V/\bT^k$ as a reflection group (cf. \cite[Prop. 1.2]{GL}), $\Gamma$ is generated by a set $\Omega$ all of whose elements can be lifted to a nice involution in $G$. Since $(G,V)$ is non-polar and indecomposable, Proposition \ref{identity_rep_indec} implies moreover that the induced representation $(\bT^k,V)$ is indecomposable as well. We now claim that $V/\bT^k$ has no boundary, so that $(\bT^k,V)$ is reduced (cf. \cite[Prop. 1.1]{GL}). Indeed, if $\partial(V/\bT^k)\neq\varnothing$, then Propositions \ref{boundary_and_copolarity_toric_actions} would yield a proper generalized section for $(\bT^k,V)$, and hence a proper generalized section for $(G,V)$, contradiction.
\edim

We now decompose the space $V$ into its isotypical components with respect to the induced representation $(\bT^k,V)$:
\beq\label{Tk_iso_comp}
V=V_1\oplus\cdots\oplus V_m,
\eeq
Observe that the induced action of $\bT^k$ on each $V_i$ cannot be trivial, since otherwise $(\bT^k,V)$ would be decomposable, and this is impossible by Lemma \ref{easy_remarks}.

\blem\label{trivial_fixed_pts}
We have that $\dim V$ is even and $\geq 2k+2$.
\elem
\bdim
What we noticed above implies that $\dim V=2m$ is even. Moreover, by Corollary \ref{k_lines} and indecomposability of $(\bT^k,V)$ (cf. Lemma \ref{easy_remarks}) we get that $V$ has $m\geq k+1$ $\bT^k$-isotypical components, hence $\dim V=2m\geq 2k+2$.
\edim

Consider now a nice involution $\omega\in G$. Since we have $0\leq\dim Z_G(\omega)\leq k$, formula (\ref{nice_involutions}) yields $1\leq\codim_VV^\omega\leq k+1$.

\blem\label{dim_ZG_k}
If $\omega\in G$ is a nice involution, then $\text{\em codim}_VV^\omega\neq 1$ or, equivalently, $\dim Z_G(\omega)\neq k$.
\elem
\bdim
Assume that $\omega\in G$ is a nice involution such that $\dim Z_G(\omega)=k$; then the identity component of $G$, $\bT^k$, is contained in $Z_G(\omega)$, and $\omega$ preserves all $\bT^k$-invariant subspaces of $V$. Since $\omega$ is a reflection, there exists a $\bT^k$-invariant subspace $U$ of $V$ such that $\omega(U)=U$ and $\omega|_U$ is a reflection. If $\dim U=2$, then $(\bT^k,U)$ is equivalent to an irreducible $S^1$-representation, which is given by rotations in $U\simeq\bR^2$; this is impossible since $Z_G(\omega)\supseteq\bT^k$. Hence $\dim U=1$, and $\bT^k$ fixes $U$ pointwisely, a contradiction.
\edim

As a corollary we immediately deduce the following:

\blem\label{no_reflections}
Let $\omega\in G$ be a nice involution. Then $0\leq\dim Z_G(\omega)\leq k-1$ and $2\leq\text{\em codim}_VV^\omega\leq k+1$.
\elem

We now assume that $\omega\in G$ is a nice involution satisfying $\codim_VV^\omega=k+1$ or, equivalently, $\dim Z_G(\omega)=0$. Then $\Ad(\omega):\gt^k\rightarrow\gt^k$ is an involution with no non-trivial fixed points and conjugation $c_\omega$ in $\bT^k$ with respect to $\omega$ coincides with the inversion map $t\mapsto t^{-1}$. If $\phi:\bT^k\rightarrow S^1$ is a homomorphism, we have then $\phi\circ c_\omega=\bar{\phi}$; this implies that $\omega$ preserves all real isotypical components of $V$ (i.e. $\omega$ preserves decomposition (\ref{Tk_iso_comp})).

\blem\label{non_triviality_on_invariant_subspaces}
If $\omega\in G$ is a nice involution such that $\text{\em codim}_VV^\omega=k+1$, and $U\subseteq V$ is a non-trivial $\bT^k$-invariant subspace of $V$ which is preserved by $\omega$, then the action of $\omega$ on $U$ is not trivial.
\elem
\bdim
If $\omega$ acts trivially on $U$ we have, for any $t\in\bT^k$ and $u\in U$,
$$
\bar{t}\cdot u=c_\omega(t)\cdot u=t\cdot u,\qquad \text{i.e.}\qquad t^2\cdot u=u. 
$$
Thus $\bT^k$ acts trivially on $U$, which implies $U=\{O\}$.
\edim

\bprop\label{codim_k+1}
If $G$ contains a nice involution $\omega$ so that $\text{\em codim}_VV^\omega=k+1$, then $V$ decomposes as
$$
V=V_1\oplus\cdots \oplus V_{k+1}
$$
where the $V_i$'s are the irreducible isotypical components of $V$; moreover $\omega$ preserves each $V_i$ and is a reflection on it. In particular we have $\dim V=2k+2$ and $\text{\em chm}(G,V)=k+2$.
\eprop
\bdim
Let $V_i$, $i=1,\dots,m$, be the $\bT^k$-isotypical components of $V$, and let $\omega_i$ be the restriction of $\omega$ to $V_i$. Since 
$$
\sum\codim_{V_i}V_i^{\omega_i}=\codim_V V^\omega=k+1,
$$
Lemma \ref{non_triviality_on_invariant_subspaces} implies $m\leq k+1$, so $m=k+1$ by Corollary \ref{k_lines}. Using again Lemma \ref{non_triviality_on_invariant_subspaces} we get $\dim V_i=2$ for all $i=1,\dots,k+1$, so $\dim V=2k+2$ and $\chm(G,V)=k+2$, as claimed.
\edim

\section{Proof of the main results}\label{Proof_of_the_main_results}
In this Section we are going to prove our main results, beginning form Theorem \ref{copolarity_1}. As mentioned in the introduction, we shall restate such a Theorem in the setting of abstract copolarity:

\bteor\label{cohomogeneity_abstract_copolarity_1}
Let $\rho:H\rightarrow\O(W)$ be a non-reduced, indecomposable representation of a connected, compact Lie group $H$ of abstract copolarity $1$. Then $\rho$ has cohomogeneity $3$.
\eteor
\bdim
Let $(G,V)$ be a minimal reduction of $\rho$; then $G$ is a finite extension of the $1$-dimensional torus $\bT^1$, and by Lemma \ref{easy_remarks} it contains a nice involution $\omega$. Applying Lemma \ref{no_reflections} we deduce $\codim_VV^\omega=2$, so $\chm(\rho)=\chm(G,V)=3$ by Proposition \ref{codim_k+1}.
\edim

Now Corollary \ref{c1_iff_ac1} easily follows:

\bdim[Proof of Corollary \ref{c1_iff_ac1}]
Since a representation has abstract copolarity $0$ if and only if it is polar (cf. \cite[p. 69]{GL}), it is clear that if $\rho:=(H,W)$ has copolarity $1$, than it has abstract copolarity $1$.

Conversely, assume that $\rho$ has abstract copolarity $1$. If it is indecomposable, then either it is reduced and there is nothing to prove, or it is not reduced and $\chm(\rho)=3$ by Theorem \ref{cohomogeneity_abstract_copolarity_1}. Now any non-polar representation of cohomogeneity $3$ has copolarity $1$ (cf. \cite{Str}), so we are done in this case.

Assume next that $\rho$ is decomposable. Exploiting Remark \ref{ac_of_a_product}, we may suppose that $\rho$ is orbit equivalent to a product representation $(H_1\times H_2,W_1\oplus W_2)$, where $(H_1,W_1)$ is polar and $(H_2,W_2)$ is indecomposable with abstract copolarity $1$. By the above discussion $(H_2,W_2)$ has copolarity $1$, so the same is true for $\rho$ since the copolarity of a product representation is the sum of the copolarities of the factors. 
\edim

We now begin the work that will lead us to the proof of Theorem \ref{cohomogeneity_abstract_copolarity_2}; in particular from now on $\rho:=(H,W)$ will denote a non-reduced, indecomposable representation of a connected, compact Lie group $H$ of abstract copolarity $2$, and $(G,V)$ a minimal reduction of $\rho$. Clearly $G$ is a finite extension of the $2$-dimensional torus $\bT^2$, and contains a nice involution $\omega$ thank to Lemma \ref{easy_remarks}. By Lemma \ref{no_reflections} we know moreover that $\codim_VV^\omega\in\{2,3\}$. If  $\codim_VV^\omega=3$, Proposition \ref{codim_k+1} implies $\chm(\rho)=\chm(G,V)=4$ and we are done; so we may assume that all nice involutions $\omega\in G$ satisfy $\codim_VV^\omega=2$.

In what follows $\omega$ will denote an arbitrarily fixed nice involution in $G$. Our previous assumption implies then that  $\Ad(\omega):\gt^2\rightarrow\gt^2$ is an involution with eigenvalues $\pm 1$. Denote by $U_\pm\subseteq\gt^2$ the $1$-dimensional eigenspace corresponding to the eigenvalue $\pm 1$, and let $\gs^\pm$ be the annihilator of $U_\pm$, which is a line in $(\gt^2)^\ast$. Finally denote by $V_\pm$ the sum of all $\bT^2$-isotypical components of $V$ which induce the line $\gs^\pm$. If $\bar{V}:=(V_+\oplus V_-)^\perp$, $V$ can be written as
\beq\label{decomposition_omega}
V=V_+\oplus V_-\oplus \bar{V}.
\eeq
Observe that decomposition (\ref{decomposition_omega}) depends on the choice of $\omega$.

\bprop\label{codim_2}
With respect to decomposition {\em (\ref{decomposition_omega})}, $V_+=\{O\}$ and $\bar{V}$ contains exactly two isotypical components $V_1$, $V_2$ which have dimension $2$. Moreover $(\bT^2,V)$ induces exactly $3$ lines in $(\gt^2)^\ast$ and $\omega$ acts as the identity on $V_-$ and interchanges $V_1$, $V_2$.
\eprop

The proof of Proposition \ref{codim_2} requires the following:

\blem\label{trivial_subrep}
Assume $W\subseteq V$ is a $\bT^2$-invariant subspace on which $\omega$ acts as the identity. Then $W\subseteq V_-$. 
\elem
\bdim
We may suppose without loss of generality that $W$ is irreducible; denoting by $\phi$ the corresponding character, we need to prove that $\d\phi_e(x)=0$ whenever $x\in U_-$. Since such an $x$ must be of the form $\Ad(\omega)(y)-y$ for a suitable $y\in\gt^2$, this corresponds to show that $\d(\phi\circ c_\omega)_e(y)=\d\phi_e(y)$ for any $y\in\gt^2$. On the other hand if $t\in\bT^2$ and $w\in W$ we have
$$
\phi(t)w=t\cdot w=\omega t\omega\cdot w=c_\omega(t)\cdot w=\phi(c_\omega(t))w,
$$
so that $\phi=\phi\circ c_\omega$ on $\bT^2$ and our claim follows.
\edim

We will denote by $\omega_\pm$, $\bar{\omega}$, respectively, the restriction of $\omega$ to $V_\pm$, $\bar{V}$.

\bdim[Proof of Proposition \ref{codim_2}]
First we note that $\bar{V}\neq\{O\}$; indeed, otherwise $(\bT^2,V)$ would induce only two lines in $(\gt^2)^\ast$, and so it would be decomposable by Corollary \ref{k_lines}, contradicting Lemma \ref{easy_remarks}. Now, let $V_1$ be a $\bT^2$-isotypical component contained in $\bar{V}$; clearly $V_2:=\omega(V_1)$ is a $\bT^2$-isotypical component contained in $\bar{V}$ different from $V_1$, so by \cite[Lemma 6.2]{GL} we deduce that the action $(\bT^2,V_i)$, $i=1,2$, has cohomogeneity $1$. Since principal orbits of such representations are $1$-dimensional, we get $\dim V_i=2$ and $\codim_{\bar{V}}\bar{V}^{\bar{\omega}}\geq\dim V_i=2$. From
\beq\label{formula_codimension}
\codim_{V_+}V_+^{\omega_+}+\codim_{V_-}V_-^{\omega_-}+\codim_{\bar{V}}\bar{V}^{\bar{\omega}}=\codim_VV^\omega=2
\eeq
we obtain then $\codim_{V_+}V_+^{\omega_+}=0=\codim_{V_-}V_-^{\omega_-}$ so $\omega$ acts as the identity both on $V_+$ and $V_-$. This implies $V_+=\{O\}$ by Lemma \ref{trivial_subrep}. Now set $U:=V_1\oplus V_2$ and denote by $\omega'$ the restriction $\omega|_U$. Since $\codim_U U^{\omega'}=2$, a splitting formula for $\codim_V V^\omega$ similar to (\ref{formula_codimension}) shows that $\omega$ acts trivially on the orthogonal complement $U^\perp$ of $U$ in $\bar{V}$, so $U^\perp=\{O\}$ by Lemma \ref{trivial_subrep} and $\bar{V}=V_1\oplus V_2$.
\edim

Let $\omega_1\in G$ be a nice involution, and decompose $V=V_-\oplus V_1\oplus V_2$ as in Proposition \ref{codim_2} with respect to $\omega_1$. We denote by $\gs^-$, $\gs_1$, $\gs_2$ the lines in $(\gt^2)^\ast$ induced by $(\bT^2,V_-)$, $(\bT^2,V_1)$, $(\bT^2,V_2)$ respectively. Any nice involution $\omega_2\in G$ permutes the $\bT^2$-isotypical components of $V$, so, via its natural action on $(\gt^2)^\ast$, it permutes the lines $\gs^-$, $\gs_1$, $\gs_2$. If $\omega_2$ does not fix $\gs^-$, then $\dim V_-=2$, $\dim V=6$ and $\chm(G,V)=4$, so Theorem \ref{cohomogeneity_abstract_copolarity_2} holds. Note that in this case all $\bT^2$-isotypical components of $V$ are irreducible and $G$ acts as the full permutation group on them, so $(G,V)$, as well as the original representation $(H,W)$, is irreducible. Henceforth we shall suppose that all nice involutions $\omega\in G$ fix $\gs^-$ and interchange $\gs_1$, $\gs_2$. 

\boss
Under the above assumption it is not hard to prove that any two nice involutions in $G$ project onto the same element in the quotient $G/G^\circ$, thus $G/G^\circ=\bZ_2$ in this case.
\eoss

Lemma \ref{trivial_fixed_pts} implies now that $\dim V_-\geq 2$; if $\dim V_-=2$, then $\dim V=6$ and $\chm(H,W)=4$, so we shall assume $\dim V_-\geq 3$ and show that we are in case (ii) of Theorem \ref{cohomogeneity_abstract_copolarity_2}.

Notice that $V_1\oplus V_2$ and $V_-$ are $G$-invariant subspaces of $V$; by \cite[Lemma 5.1]{GL}, we can then find two $H$-invariant orthogonal complementary subspaces $W_1$, $W_2\subseteq W$ such that 
$$
W_1/H=(V_1\oplus V_2)/G,\quad W_2/H=V_-/G.
$$
Let $\rho_i$ be the induced representaton $(H,W_i)$; recalling that $\rho:=(H,W)$ in the notation of Theorem \ref{cohomogeneity_abstract_copolarity_2}, we get $\rho=\rho_1\oplus\rho_2$.

First observe that $(G,V_-)$ has no non-trivial fixed points by Lemma \ref{trivial_fixed_pts}; so $\rho_2$ has no non-trivial fixed points either (cf. \cite[Section 5]{GL}). Moreover, since any nice involution $\omega\in G$ acts trivially on $V_-$, $(G,V_-)$ has a $1$-dimensional kernel (the line $\gs^-$), and so it has the same orbits as $(S^1,V_-)$. Since $\dim V_-\geq 3$, $(S^1,V_-)$ cannot be polar and is a minimal reduction of $(H,W_2)$. Moreover, applying Proposition \ref{boundary_and_copolarity_toric_actions}, we deduce that $V_-/S^1=W_2/H$ does not have boundary. This implies, by \cite[Prop. 1.1]{GL}, that $\rho_2$ is itself reduced, so $\rho_2(H)=S^1\subseteq\SO(W_2)$ is $1$-dimensional.

We now study $(G,V_1\oplus V_2)$. Clearly it is irreducible, has cohomogeneity $2$ and so is polar. This is of course true for $\rho_1$ as well, so the latter is orbit-equivalent to the isotropy representation of a rank $2$ symmetric space (cf. \cite{Dad}). The proof that this is in fact a real grassmannian will require a few remarks.

We begin with the observation that $\rho_2(H)=S^1$, hence $H$ cannot be semisimple; we shall write then $H=Z^\circ\cdot H_s$, where $Z^\circ$, the identity component of the centre $Z(H)$ of $H$, is a torus $\bT^a$, $a\geq 1$, and $H_s$ is the semisimple part of $H$. The same argument implies moreover that the induced representation $(H_s,W_2)$ is trivial.

\blem\label{Lemma1}
$\rho_1(H)$ has $1$-dimensional centre.
\elem
\bdim
Clearly $\rho_1(Z^\circ)\subseteq Z(\rho_1(H))^\circ$, and the latter is either $0$ or $1$-dimensional. Assume by contradiction that $\dim Z(\rho_1(H))=0$. Then $Z^\circ\subseteq\ker\rho_1$, which implies that $(H,W)$ is decomposable; indeed in this case $(H,W)$ is orbit-equivalent to the product representation $(H_s\times Z^\circ,W_1\oplus W_2)$. 
\edim

\blem\label{Lemma2}
Set $K:=\rho_1(H)$, and write $K=S^1\cdot K_s$, where $K_s$ denotes the semisimple part of $K$. Then $(K,W_1)$ is not orbit-equivalent to the induced representation $(K_s,W_1)$.
\elem
\bdim
Representations $(K,W_1)$ and $\rho_1=(H,W_1)$ are orbit-equivalent, and so are $(K_s,W_1)$ and $\rho_s:=(H_s,W_1)$ since $K_s=\rho_1(H_s)$. It is then enough to prove that $\rho_1$ and $\rho_s$ are not orbit-equivalent. This is easily done, since it is not hard to show that if this were not the case, then $(H,W)$ and the product representation $(H_s\times Z^\circ,W_1\oplus W_2)$ would have the same orbits.
\edim

Lemma \ref{Lemma1} and the main theorem in \cite{EH} imply that $\rho_1$ is orbit-equivalent to the isotropy representation of a rank $2$ irreducible Hermitian symmetric space, i.e. one of the following (cf. \cite[p.518-520]{Hel}):
\bit
\item AIII($p\geq 3$, $q=2$): $\SU(p+2)/\text{S}(\U(p)\times\U(q))$,
\item DIII($n=5$): $\SO(10)/\U(5)$,
\item BDI($p\geq 3,\; q=2$): $\SO(p+2)/(\SO(p)\times\SO(2))$,
\item EIII: $\text{\bf E}_6/\SO(10)\cdot \bT^1$.
\eit
Only the isotropy representation of BDI($p\geq 3,\; q=2$) satisfies condition in Lemma \ref{Lemma2} (cf. \cite{EH}), and the first part of Theorem \ref{cohomogeneity_abstract_copolarity_2} is finally proved.

\boss
Assuming $\rho$ faithful, we have $a=\dim Z(H)=1$. Indeed, since $\rho_1(H)=K=K_s\cdot S^1$ and $\rho_2(H)=S^1$, if $a\geq 3$ then $Z^\circ\simeq\bT^a$ contains a subgroup of positive dimension acting trivially on $W$. So $a\leq 2$. Suppose now $a=2$. Since $\rho_1(Z^\circ)\simeq\rho_2(Z^\circ)\simeq S^1$ and $\rho$ is faithful, we may write $Z^\circ=S^1_1\times S^1_2$ for suitable subgroups $S^1_i$ of $Z^\circ$ isomorphic to $S^1$ in such a way that $S^1_1$ acts trivially on $W_2$ and $S^1_2$ acts trivially on $W_1$. Hence $H_s\cdot S^1_1$ acts trivially on $W_2$ and we conclude that $\rho$ is orbit-equivalent to the natural product representation $((H_s\times S^1_1)\times S^1_2, W_1\oplus W_2)$. Once again this is a contradiction.
\eoss

We now turn to the proof of the second part of Theorem \ref{cohomogeneity_abstract_copolarity_2}; in what follows we shall use notation introduced in its statement. We set moreover $H:=\SO(2)\times\SO(n)$, $W_1:=\bR^2\otimes\bR^n$, and we define $W_2$ to be the representation space of $\rho_2$. In this way we can write
$$
\rho=(H,W),\qquad \rho_1=(H,W_1),\qquad\rho_2=(H,W_2),
$$
where $W:=W_1\oplus W_2$.  Here are some preliminary results:

\blem\label{Lemma3}
Copolarity and abstract copolarity of $\rho$ are both $\leq 2$.
\elem
\bdim
The proof consists in computing the Luna-Richardson-Straume reduction of $\rho=\rho_1\oplus\rho_2$. It is easily seen that the principal isotropy group $K_1$ of $\rho_1$ is isomorphic to $\bZ_2\times \SO(n-2)$, so the principal isotropy group $K$ of $\rho$ is isomorphic either to $\SO(n-2)$ or to $\bZ_2\times \SO(n-2)$. In both cases $N_H(K)^\circ\simeq\bT^2\times\SO(n-2)$ and  $\left(\!\frac{N_H(K)}{K}\!\right)^{\!\circ}\simeq\bT^2$, as claimed.
\edim

\blem\label{Lemma4}
$\rho$ is indecomposable. 
\elem
\bdim
Suppose this is not true; we may assume without loss of generality that $\rho$ is orbit-equivalent to
$$
\tilde{\rho}:=(H'\times H'', W'\oplus W''),
$$
where $W_1\oplus W_2=W'\oplus W''$ and $W'$, $W''$ are $G$-invariant with positive dimension. Since $W_1$ is an irreducible isotypical component of $\rho$, we may assume also that $W_1\subseteq W'$, and that
$$
W'=W_1\oplus\bar{W}_1.
$$
Here $\bar{W}_1$ is a sum of irreducible components of $\rho$ and is $H'$-invariant, while $\bar{W}_1\oplus W''=W_2$ is a sum of isotypical components of $\rho$. Using that $\rho$ and $\tilde{\rho}$ have the same orbits, we easily conclude that $(H,W_2)$ and the product representation $(H'\times H'',\bar{W}_1\oplus W'')$ are orbit-equivalent; since principal orbits of $(H,W_2)$ have dimension $1$, this forces one of the two representations $(H',\bar{W}_1)$, $(H'',W'')$ to be trivial. Since $\dim W''\geq 1$ and $(H,W_2)$ has no non-trivial fixed points, this implies $\bar{W}_1=\{O\}$, and thus $(H',W')$, $(H'',W'')$ are orbit-equivalent to $\rho_1$, $\rho_2$ respectively. It follows that
$$
\rho\simeq_{o.e.} (H'\times H'', W'\oplus W'')\simeq_{o.e.} 
  ((\SO(2)\times\SO(n))\times S^1,(\bR^2\otimes\bR^n)\oplus W_2)
$$
which is absurd since principal orbits of $((\SO(2)\times\SO(n))\times S^1,(\bR^2\otimes\bR^n)\oplus W_2)$ and of $\rho$ don't have the same dimension.
\edim

\blem\label{chm_5}
We have $\text{\em chm}(\rho)\geq 6$.
\elem
\bdim
Clearly (cf. \cite[p.92]{GL}) $\chm(\rho)=2+\dim W_2$. Since $\rho_2$ is not polar and has no non-trivial fixed points, $\dim W_2\geq 4$ and we are done.
\edim

We may finally prove the second part of Theorem \ref{cohomogeneity_abstract_copolarity_2}:

\bdim 
Using Lemmas \ref{Lemma4} and \ref{chm_5} it is enough to prove that c$(\rho)=\ac(\rho)=2$.

First we observe that $\rho$ is not polar (since otherwise both $\rho_1$ and $\rho_2$ would be polar, \cite{Dad}), and its abstract copolarity cannot be $1$ (otherwise we would have $\chm(\rho)=3$ by Theorem \ref{cohomogeneity_abstract_copolarity_1} and Lemma \ref{Lemma4}). Hence, the abstract copolarity of $\rho$ has to be at least $2$, which, together with Lemma \ref{Lemma3}, implies that the abstract copolarity (and of course the copolarity) of $\rho$ is exactly $2$.
\edim

\begin{small}

$\,$

$\,$

Francesco Panelli 

{\scshape Dipartimento di Matematica e Informatica "Ulisse Dini",}

{\scshape Universit\`a di Firenze,}

{\scshape Viale Morgagni 67/A, 50134 Firenze, Italy}

{\em E-mail address}: \texttt{francesco.panelli@unifi.it}

{\em E-mail address}: \texttt{francescopanelli@virgilio.it}

$\,$

$\,$

Carolin Pomrehn

{\em E-mail address}: \texttt{cpomrehn@googlemail.com}

\end{small}
\end{document}